\documentclass[twoside,leqno]{article}

\usepackage[letterpaper]{geometry}
\usepackage{ltexpprt}
\usepackage{hyperref}

\usepackage{soul}  
\usepackage{amsfonts}
\usepackage{graphicx}

\begin{document}

\newcommand{\liz}[1]{{\color{red}\textbf{Liz says:} #1}}
\graphicspath{{Figures/}}
\renewcommand{\epsilon}{\varepsilon}
\newcommand{\e}{\varepsilon}

\newcommand{\ismail}[1]{{\color{blue}\textbf{Ismail answer:} #1}}
\newcommand{\updated}[1]{{\color{purple}\textbf{ #1 }}}
\newtheorem{definition}{Definition}[section]

%


\title{\Large A Case Study on Identifying Bifurcation and Chaos with CROCKER Plots}
\author{
\.{I}smail G\"{u}zel
\thanks{Department of Mathematics, \.{I}stanbul Technical University, Maslak, \.{I}stanbul, Türkiye (\textit{iguzel@itu.edu.tr}).}
\and
Elizabeth Munch
\thanks{Dept of Computational Mathematics, Science and Engineering; Dept of Mathematics, Michigan State University, East Lansing, MI (\textit{muncheli@msu.edu}).}
\and
Firas Khasawneh
\thanks{Mechanical Engineering, Michigan State University, East Lansing, MI (\textit{khasawn3@egr.msu.edu}).}
}
\date{}

\maketitle


\fancyfoot[R]{\scriptsize{Copyright \textcopyright\ 2022 by SIAM\\
Unauthorized reproduction of this article is prohibited}}





\begin{abstract} \small\baselineskip=9pt
The CROCKER plot is a coarsened but easy to visualize representation of the data in a one-parameter varying family of persistence barcodes. 
In this paper, we use the CROCKER plot to view changes in the persistence under a varying bifurcation parameter. 
We perform experiments to support our methods using the Rössler and Lorenz system and show the relationship with common methods for bifurcation analysis such as the Lyapunov exponent.
\end{abstract}

\section{Introduction}

    Dynamical systems having chaotic behavior can be found in a variety of domains, including biology, economics, and engineering. 
    Of great interest is classifying regimes of a dynamical system's parameter space based on whether the resulting output is chaotic or not. 
    Of course, many options are available for such a question; in this work we focus on a data driven approach using topological data analysis (TDA) to measure the shape and structure of the attractor of the system as a proxy for behavior. 
    The idea of combining methods from TDA with dynamical systems and/or time series analysis is not new 
	\cite{khasawneh2017utilizing,myers2019persistent,charo2020topologyLagrangian,tymochko2020using,gradingmodels2021,maletic2016persistent,TACTS,mittal2017topological}. 
    In this work, we focus on a new way to encode information about the bifurcation, namely the CROCKER plot \cite{crocker-Lori2015,crocker-Lori2019,crocker-stack}.
    This construction is a simplified discretization of the information from a family of persistence barcodes; but which, as we will see here, still has a great deal of information related to the behavior of the system. 
    We further show that the even further simplified information encoded in the $L_1$ norm of the Betti vectors used to construct the CROCKER plots is itself related to the bifurcation information via the Lyapunov exponent. 
    
\section{Materials and Methods}
	
    \subsection{The Lyapunov exponent and chaos}
    
    One way to measure whether the solution to a given dymaical system is chaotic or not is to compute the maximum Lyapunov exponent $\lambda$.
    This value is the mean rate of exponential divergence or convergence of two neighboring initial points in the phase space of a dynamical system. 
    It is a quantitative measure where $\lambda>0$ implies the system is chaotic, $\lambda = 0 $ implies periodic, and $\lambda <0$ implies stable. 
    This value can be approximated numerically but still requires considerable computational effort.
    
\subsection{Persistence Barcodes}
    Topological Data Analysis (TDA) is a collection of tools that have the ability to measure the shape of given data. 
    One of the popular tools of TDA is persistent homology which stores information about shape by encoding how a parameterized space's homology changes over a varying parameter. 
    The structure measured depends on the dimension of homology used: connected components are encoded in zero-dimensional homology, loops are in one-dimensional homology, and voids in two-dimensional homology. 
    We briefly introduce the relevant background here, but invite the interested reader to seek out  a more thorough treatment \cite{Dey2021,Oudot2017a}.
    
    To compute persistence on a point cloud $X$ with distance $d$, such as those taken as input from a dynamical system, we need to construct a simplicial complex. One such option for this task is the Vietoris-Rips complex, defined as follows.

	\begin{definition}
	    Given a point cloud $ X $, the \emph{Vietoris-Rips complex} is  the simplicial complex  whose simplices are built on vertices that are at most $\e$  apart,
		\[ R_{\varepsilon}(X) = \{\sigma \subset X\mid d(x,y)\leq \varepsilon, \mbox{ for all }x,y\in \sigma\}. \]
	\end{definition}
    
	Note that this construction is dependent on the choice of a proximity parameter $\e$. 
	If we choose an increasing sequence of proximity parameters, $ \epsilon_0 \leq \epsilon_1\leq \cdots \leq \epsilon_n $, we obtain a filtration as	$ R_{\epsilon_0} \subseteq R_{\epsilon_1} \subseteq \cdots \subseteq R_{\epsilon_n}.$ 
	We then compute the $p$-dimensional homology of each complex, which is a vector space $H_p(R_{\e_i})$. 
    Fundamental results in algebraic topology mean that we also get  linear maps between each of them; that is,
	$ H_{p}(R_{\epsilon_0}) \rightarrow H_{p}(R_{\epsilon_1}) \rightarrow \cdots \rightarrow H_{p}(R_{\epsilon_n})$. 
	This is known as a persistence module. 
	
\begin{figure}
    	\centering
    	\includegraphics[scale=0.4]{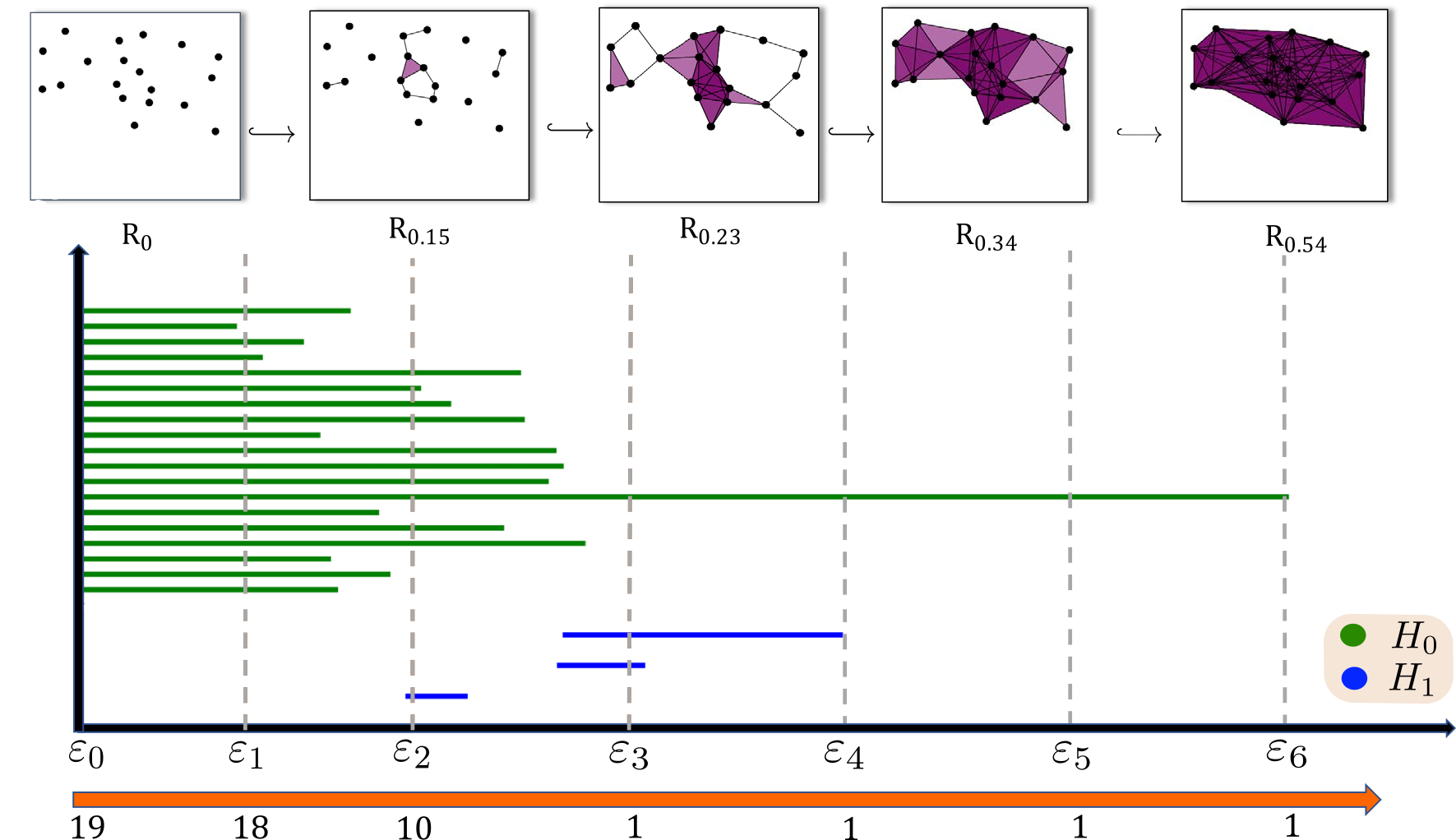}
    	\caption{The Vietoris Rips filtration and Betti numbers corresponding to given filtration parameter $\varepsilon$ on the orange line. One can follow the above either two ways to get Betti vector for the dimension $p=0$  with the partition $P=\{0,0.7,0.15,0.23,0.34,0.44,0.54\}$, $ Bv_0(X;P)=(19,18,10,1,1,1,1).$}\label{bettibarcode}
    \end{figure}

    The information in a persistence module can be uniquely represented through a collection of pairs $(\e_{birth},\e_{death})$ in a decomposition of the module, where each pair represents the parameter values for which a homological feature appeared and disappeared. 
    We visualize this information in a persistence barcode, i.e.~a collection of horizontal line segments as in the example of Fig.~\ref{bettibarcode}. 
    We place bars on the vertical axis (where order does not matter) while the horizontal axis represents the life span of each homology class in terms of the parameter $ \varepsilon $. 
    For the purposes of our work, however, we will not need the full barcode but instead focus on a simplified invariant, the Betti vector, which define here in the special case of the data we are using.

    
	\begin{definition}{\label{defn:Bettinumber}}
      For a fixed $\e$ and $p$, the $ p^{th} $
      \emph{Betti number} $ \beta_p (R_{\epsilon}) $ is 
      the dimension of the $ p^{th} $ homology group of $ R_{\epsilon} $, i.e.
      \[ \beta_p(R_{\epsilon}) = \dim H_p(R_{\epsilon}). \] 

      Let $P=\{\epsilon_0, \epsilon_1,\dots,\epsilon_{N}\}$ be a partition of the interval $(0, \epsilon_{N})$. 
      The $p^{th}$ dimensional \emph{Betti vector} is the ordered sequence of the $p^{th}$ dimensional Betti numbers, that is
        $$
        Bv_p(X;P) = (\beta_p(R_{\e_0}),\beta_p(R_{\e_1}),\dots, \beta_p(R_{\e_N})).
        $$
    \end{definition}

    In the example of Fig.~\ref{bettibarcode}, when we draw a vertical line at a particular $ \varepsilon $, the number of intersecting bars of a particular homological dimension is the rank of the corresponding homology group, i.e.~the Betti number for that parameter $\e$. 
    We can see that for a fixed $\e$, $ \beta_{0} $ counts the number of connected components and $ \beta_{1} $, counts the number of loops.

\subsection{Time-Evolving Metric Spaces}
    In the previous section, we defined a varying parameter given fixed input data to compute a Betti vector.  
    We next consider the case of an evolving series of point clouds. 
    While several methods for handling this sort of input data with persistence have been developed \cite{vineyards,multiparameterrank,move-schedules},
    in this work we focus on the CROCKER plot \cite{crocker-Lori2015,crocker-Lori2019,crocker-ML2019,crocker-stack}.
        
	\begin{definition}
	    For a given collection of point clouds $ \mathcal{X} = \{X_1, X_2, \cdots, X_T\} $, the \emph{CROCKER matrix} of this collection  is given by
        \[ 
        \mathrm{CR}(\mathcal{\mathcal{X}}) = ( Bv(X_1;P), Bv(X_2;P), \dots, Bv(X_T;P) ),
        \]
        where $Bv(\bullet)$ is the $p^{th}$ dimensional Betti vector for the partition $P = \{\epsilon_0,\epsilon_1,\dots,\epsilon_N\}$.
	\end{definition} 
    
    Since the CROCKER matrix is a two dimensional array, we can visualize it in terms of a heat map, often called a \textit{CROCKER plot}; see Fig.~\ref{RosslerResult} for an example.

\subsection{Algorithm}

    In this paper, we investigate the CROCKER plots resulting from a dynamical system changing over an input parameter, and compare them to information from the bifurcation diagram and Lyapunov exponent information. 
    Our algorithm is the following, where we repeat the following process for each investigated parameter on the dynamical system:
    \begin{itemize}
        \item Obtain the states of the nonlinear system with python library \texttt{teaspoon}~\cite{teaspoon}.
        \item Calculate the barcode using the python library \texttt{ripser}~\cite{ripser} with greedy sub-sampling algorithm~\cite{greedyalgo}.
        \item Find the maximum death time $d_{max}$ for each homological dimension $p\in\{0,1\}$.
        \item Get 100 equally-spaced values of the proximity parameter $\varepsilon$ between 0 and $d_{max}$.
        \item Obtain Betti vectors for each $p\in\{0,1\}$.
    \end{itemize}

\section{Results}
    We show the pipeline as applied to a nonlinear dynamical system known as the Rössler system, given by 
    \begin{equation}
    	\dot{x} =  -y -z,\qquad
		\dot{y} =  x + ay,\qquad
		\dot{z} =  b+ z(x-c)    
    \end{equation}
    In this work, we fix parameters $b=2$ and $c=4$; then vary the control parameter $a$ with the initial conditions $[x_0, y_0, z_0] = [-0.4,0.6,1]$.
	
\begin{figure}
		\centering
		\includegraphics[scale=0.27]{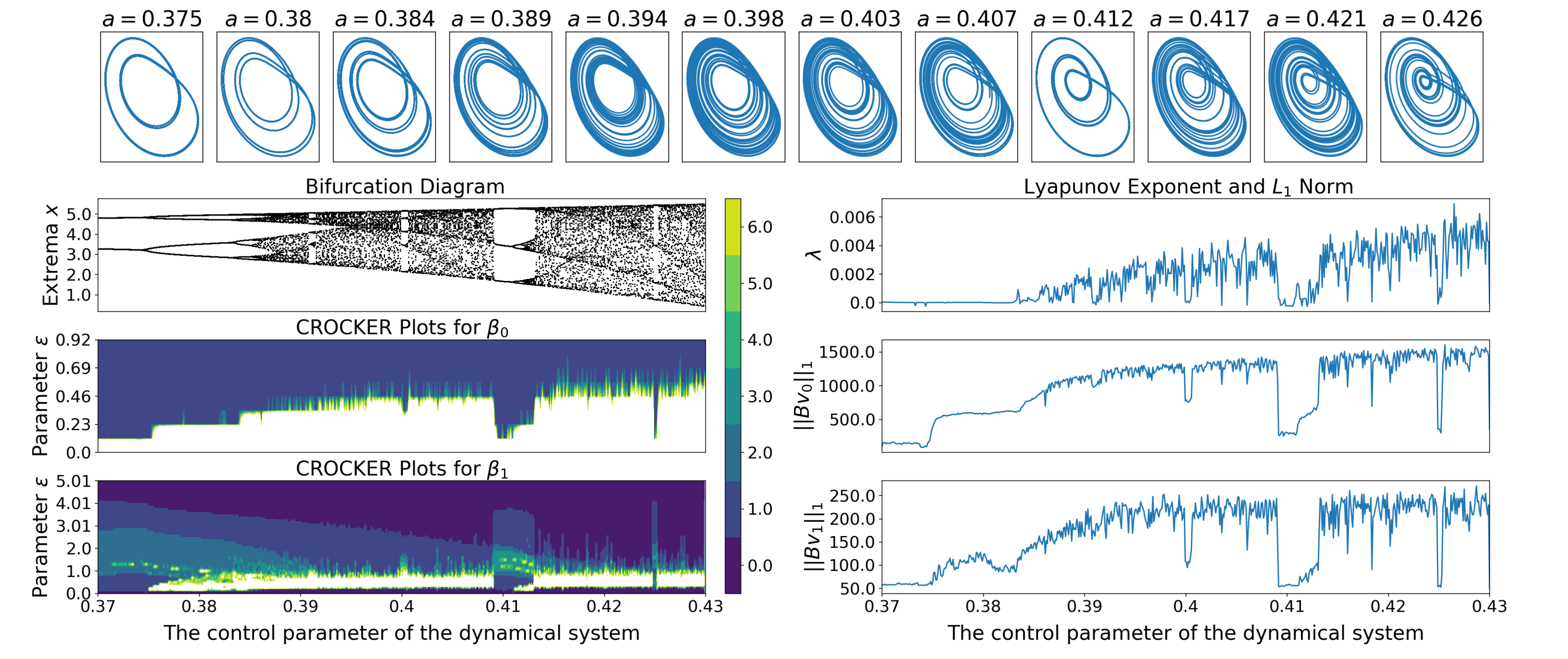}
		\caption{The CROCKER plots (left) and $L_1$ norm of each Betti vector (right) corresponding to varying the control parameter $a$ for the Rössler system.}\label{RosslerResult}
	\end{figure}

	In Fig.~\ref{RosslerResult} (left), we show the bifurcation diagram for varying parameter $a$ with 600 equally-spaced values between 0.37 and 0.43 with the other parameters fixed. 
	Below this panel, we show the CROCKER plots for homological dimensions 0 and 1 for comparison to the bifurcation diagram. 
	In the CROCKER plots, one can read the Betti numbers by using color bar which is near the figures. 
	For instance, if we consider the parameter $a$ around 0.41, there is a distinct change in structure for both the 0- and $1$-dimensional CROCKER plots.
	
	Despite the fact that the CROCKER plot is a simplified version of viewing the data than the persistence barcodes setup, we can actually choose to look at our data in an even more simplified fashion. 
    We look at  the correlation between Lyapunov exponent and the $L_1$ norm of a given Betti vector for partition $P = \{\e_0,\cdots,\e_N\}$, 
    \begin{equation}
        \| Bv_p(X;P)\|_1 = \sum \beta_p(R_{\e_i}).
    \end{equation}
    In this case, our CROCKER vector information is simplified down to a single number, reducing the CROCKER plot down to $\mathbb{R}$-valued function parameterized by $a$. In Fig.~\ref{RosslerResult} (right), we demonstrate the Lyapunov exponent and the $L_1$ norms of the Betti vectors. The computed Pearson coefficient values between the Lyapunov exponent and the $L_1$ norms for 0- and 1-dimensional information was 0.86 and 0.83, respectively.

    We also demonstrate our methods on the  Lorenz system which consists of three ordinary differential equations referred to as Lorenz equations:
    \begin{equation}
        \dot{x} = \sigma(y-x),\quad
        \dot{y} = x(\rho-z)-y,\quad
        \dot{z} = xy-\beta z.
    \end{equation}
    The fixed parameters $\sigma = 10,\: \beta = 8/3$ and the control parameter $\rho$ varying on 600 equally-spaced values between 90 and 105 with the initial conditions $[x_0, y_0, z_0] = [10^{-10},0,1]$. 
    
\begin{figure}[ht]
		\centering
		\includegraphics[scale=0.27]{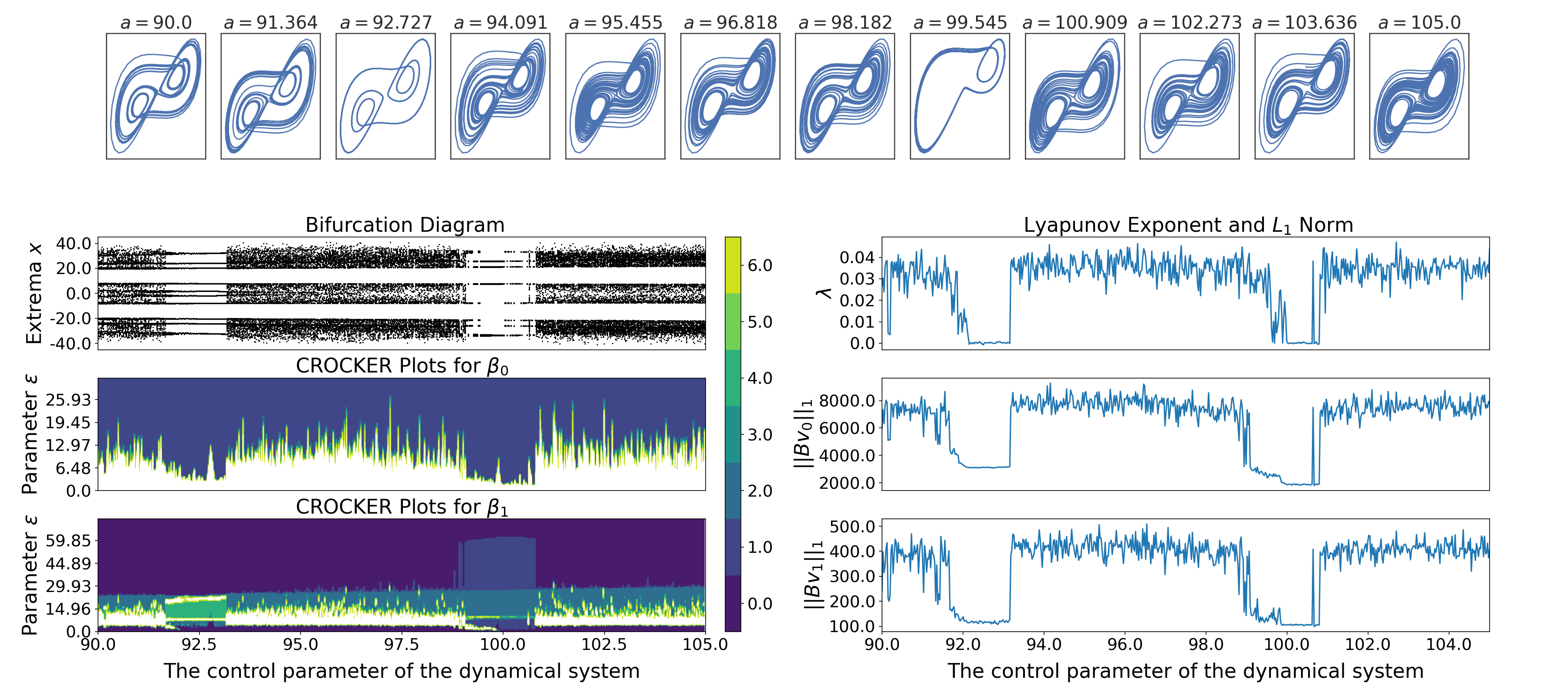}
		\caption{The CROCKER plots (left) and $L_1$ norm of each Betti vector (right) corresponding to varying the control parameter $\rho$ for the Lorenz system.}\label{LorenzResult}
	\end{figure}

    In Fig. \ref{LorenzResult} (left), as in the Rössler system, Lorenz's CROCKER and $L_1$ norms exhibit similar characteristics to the bifurcation diagram and the Lyapunov exponent. However, when we examine two different parameters in the same system, some differences emerge in the case of 1-dimensional CROCKER. In particular, consider the parameters $\rho$ around 92.5 and 100. While in the 0-dimensional CROCKER, there is not an obvious difference  between the two regions, the 1-dimensional CROCKER  shows a stark contrast. For example, there are 4 noticeably persistent points around $\rho=92.5$, and there is an extremely long lived persistence bar around $\rho=100$.

    We also note that there is a clear relationship between the Lyapunov exponent and the $L_1$ norm of the CROCKER vectors as seen on the right of Fig. \ref{LorenzResult}.
    The computed Pearson coefficient values between the Lyapunov exponent, and the $L_1$ norms for 0- and 1-dimensional vectors was 0.85 in both cases.

    We further note the computational cost of each method.  All the computations we performed on a Ubuntu 20.10 desktop with 16 GB RAM, Intel(R) Core(TM) i7-9700 CPU 3.00GHz, and 8 cores using the python language. We use the example of the Rössler system with the parameters $[a,b,c] = [0.41,2,4]$ and the initial conditions $[-0.4,0.6,1]$. To compute the Lyapunov exponent, it took $716 \pm 17.3$ ms  over 7 runs while computing the $L_1$ norm of Betti vector took $562 \pm 12.3$ ms and $566 \pm 14.7$ ms  over 7 runs for the dimension 0 and 1, respectively. So, in this case TDA is computationally less expensive than the computing the Lyapunov exponent.

\section{Discussions}
    In this work, we have begun an investigation of the use of CROCKER plots for bifurcation analysis in dynamical systems. 
    We show that in a simple test case, there is clearly a relationship between a representation of behavior in the system (the Lyapunov exponent) with the structure of the CROCKER plot, as well as with the $L_1$ norm of each Betti vector.
    
    This work, of course, leads to many interesting open questions. 
    For starters, more work must to be done to understand how changes in the CROCKER plot show up in different dynamical systems.
    Rössler and Lorenz are an easy test cases to start with due to their circular structure, so it will be interesting to see what sorts of changes can be found in the CROCKER plot given different changes in structure. 
    
    Second, as with all methods using persistent homology, there is a computational cost involved in using methods from TDA. While we have seen improvements in speed over Lyapunov, at least in part, this is because to our knowledge all code available computes Betti vectors by first computing the full persistence barcodes. 
    Might there be a more direct computation method which provides speedups relative to the desired refinement of the partition?
    
    The positives gained in using persistence are those of stability. 
    In particular, Lyapunov is notorious for its slow computation and its sensitivity to error and noise. 
    On the other hand, persistence comes with a theoretically grounded framework that should extend to the CROCKER framework in this case. Even though the CROCKER plots are unstable because of the instability of Betti numbers~\cite{crocker-stack, johnson2021instability}, the $L_1$ norm of further refined partitions of the $\epsilon$ parameter might be able to mitigate the damage. This is beyond the scope of this paper, and will be an interesting direction for future work.

\section*{Funding}
    The research of the first author was supported by a grant program (B\.{I}DEB 2214-A:1059B142000135) from T\"{U}B\.{I}TAK, Scientific and Technological Research Council of Turkey. This material is based upon work supported by the Air Force Office of Scientific Research under award number FA9550-22-1-0007.
\section*{Acknowledgments}  
    We also thank the three anonymous reviewers whose comments helped improve the quality of this study.

\bibliographystyle{siamplain}
\bibliography{references}
\end{document}